\documentclass{beamer}

\usepackage{ngerman}
\usepackage{bibgerm} 
\usepackage{enumerate}
\usepackage{color}

\usepackage{amsmath}
\usepackage{graphicx}
\usepackage[arrow, matrix, curve]{xy}
\usepackage{longtable}

\newtheorem{beispiel}{Beispiel}
\newtheorem{korollar}{Korollar}

\newcommand{\N}{\mathbb{N}}
\newcommand{\Z}{\mathbb{Z}}

\renewcommand {\Re}{{\rm I \kern-2pt R}}
\newcommand {\Nat}{{\rm I \kern-2pt N}}
\title{Das Kapitulationsproblem in der Klassenk\"orpertheorie}
\author{Tobias Bembom}
\date{2 April 2012}
\institute{Universit\"at G\"ottingen}
\begin{document}
\maketitle

\begin{frame}[shrink]
\begin{center}
\Large \textbf{Themen\"ubersicht}
\end{center}
\pause
\begin{itemize}[<+->]
\item Grundlagen der Klassenk\"orpertheorie
\vspace{3mm}\item Das Kapitulationsproblem \& Einordnung in die Literatur
\vspace{3mm}\item Chevalleys Theorem \& Implikationen
\vspace{3mm}\item Wachstum von Idealklassen
\vspace{3mm}\item Struktur des Kapitulationskerns \& $G$-Wirkung
\vspace{3mm}\item Kapitulation in imagin\"ar quadratischen K\"orpern \& Heuristik
\end{itemize}

\end{frame}

\begin{frame}[shrink]
\begin{center}
\large \textbf{Hauptsatz der Klassenk\"orpertheorie}
\end{center}
\pause
\begin{theorem} Sei $K$ ein algebraischer Zahlk\"orper mit Idealklassengruppe $Cl(K)$. Dann existiert eine eindeutige K\"orpererweiterung $H(K)/K$, so dass:
\vspace{2mm} \noindent \\ (i) $H(K)$ ist die maximale unverzweigte abelsche Erweiterung von $K$;
\\ (ii) $Cl(K)\cong Gal(H(K)/K)$.
\vspace{4mm}
\pause
\begin{itemize}[<+->]
\item Man nennt $H(K)$ den {\bf Hilbertklassenk\"orper} von $K$. 
\item Der Isomorphismus in (ii) wird vom {\bf Artin Symbol} induziert. 
\end{itemize}
\end{theorem}
\end{frame}

\begin{frame}[shrink]
\begin{center}
\large \textbf{Das Kapitulationsproblem und sein Historischer Hintergrund}
\end{center}
\pause
Sei $L/K$ eine Erweiterung algebraischer Zahlk\"orper.
\pause
\begin{itemize}[<+->]
\item Dann definieren wir den {\bf Lift} von Idealen von $K$ zu $L$ wie folgt:
\[\imath_{L/K}:\ \mathfrak{J}_K\rightarrow \mathfrak{J}_L,\ I\mapsto I\cdot \mathcal{O}_L\]
\vspace{2mm}
\item Dies ist ein injektiver Gruppenhomomorphismus, der kanonisch den Lift von Idealklassen induziert:
\[\bar{\imath}_{L/K}:\ Cl(K)\rightarrow Cl(L),\ [I]\mapsto [I\cdot \mathcal{O}_L],\]
wobei $I$ ein Ideal in $K$ ist und $[I]$ die von $I$ erzeugte Idealklasse. 
\vspace{4mm} \item $\bar{\imath}_{L/K}$ ist i.A. nicht mehr injektiv, d.h. es existieren Ideale in $K$, die in $K$ nicht Hauptideal sind, aber beim Liften zu $L$ zu Hauptidealen werden.
\end{itemize}
\end{frame}

\begin{frame}[shrink]
\begin{definition} Man sagt, dass solche Ideale wie oben in $L/K$ {\bf kapitulieren}. Entsprechend definieren wir den {\bf Kapitulationskern} von $L/K$ als
\[P_K(L)=ker(\bar{\imath}_{L/K}:\ Cl(K)\rightarrow Cl(L)).\]
\end{definition}
\pause
\begin{theorem}[Hilbert 94,1897]
\label{ubounds} 

Sei $L/K$ eine zyklische unverzweigte Erweiterung mit Galoisgruppe $G$. Dann gilt
\[|P_K(L)| = [L:K]\cdot |\mathcal{H}^0(G,\mathcal{O}_L^{*})|,\]
wobei $\mathcal{O}_L^{*}$ die Einheitengruppe vom Ganzheitsring von $L$ ist und $\mathcal{H}^0(G,\mathcal{O}_L^{*})=\mathcal{O}_K^{*}/N_{L/K}(\mathcal{O}_L^{*})$ die 0-te Kohomologiegruppe von $(G,\mathcal{O}_L^{*})$ ist. 
\end{theorem}
\pause
\begin{beispiel} Falls $K$ imagin\"ar quadratisch und $[L:K]=p^n$, $p>2$, dann: $|P_K(L)|=[L:K]$. 
\end{beispiel}
\end{frame}

\begin{frame}[shrink]
\vspace{4mm}
\begin{theorem}[Suzuki,1991] Sei $L/K$ eine unverzweigte abelsche Erweiterung. Dann gilt:
\[[L:K]\ \mbox{teilt}\ |P_K(L)|.\]
\end{theorem}
\pause
\begin{itemize} [<+->] 
\item Als Spezialfall hatte Furtw\"angler 1932 den {\bf Hauptidealsatz} bewiesen: Alle Ideale in $K$ werden in $H(K)$ zu Hauptidealen. 
\vspace{4mm} \item Entscheidender Baustein f\"ur den Satz von Suzuki ist die sogenannte {\bf Verlagerungstheorie}:
\end{itemize}
\end{frame}

\begin{frame}[shrink]
\begin{theorem}[Artin, 1930] Sei $L/K$ eine unverzweigte abelsche Erweiterung. Dann kommutiert das folgende Diagramm:
\[
\begin{xy}
\xymatrix{Cl(L) \ar[r]     &   Gal(H(L)/L)  \\
      Cl(K) \ar[r]\ar[u]_{\imath_{L/K}}         &   Gal(H(K)/K)\ar[u]_{Ver_{L/K}}}
\end{xy}
\]
\end{theorem}
\pause
\begin{itemize}[<+->]
\item Die horizontalen Abbildungen sind durch die Artinsymbole von $H(K)/K$, bzw. $H(L)/L$ induziert.
\vspace{3mm}\item $Ver_{L/K}$ ist die sogenannte Gruppenverlagerung von $L/K$.
\vspace{3mm}\item Durch den Satz von Artin wird das Kapitulationsproblem auf ein rein gruppentheoretisches Problem zur\"uckgef\"uhrt.
\end{itemize}
\end{frame}

\begin{frame}[shrink]
\begin{center}
\Large {\bf Interessante Fragen}
\end{center}
\pause
\begin{itemize}[<+->]
\item Was ist die genaue Ordnung des Kapitulationskerns?
\vspace{4mm} \item Suzukis Theorem macht keine Aussage \"uber die Struktur des Kapitulationskerns. Unter welchen Umst\"anden k\"onnen wir z.B. $Gal(L/K)$ in den Kapitulationskern einbetten? 
\vspace{4mm} \item Sei $K$ ein Zahlk\"orper mit beispielsweise $Cl(K)\cong C_p\times C_p$. Dann existieren $p+1$ unverzweigte zyklische Erweiterungen $L_1,...,L_{p+1}$ vom Grad $p$ \"uber $K$. Sind die Kapitulationskerne $P_K(L_i)$ irgendwie korreliert? 
\end{itemize}
\end{frame}

\begin{frame}[shrink]
\begin{center}
\Large {\bf Einschr\"ankung auf $p$-Teile} 
\end{center}
\pause
\begin{theorem} Sei $L/K$ eine Erweiterung algebraischer Zahlk\"orper. Dann gilt: 
\[N_{L/K}(\imath_{L/K}(a))=a^{[L:K]},\ \forall\ a\in{Cl(K)}.\]
\end{theorem}
\pause
\begin{itemize} [<+->]
\item Falls $([L:K],ord(a))=1$, dann $ord(a)=ord(\imath_{L/K}(a))$.
\vspace{4mm} \item K\"onnen uns auf den $p$-Teil der Klassengruppe von $K$ und auf unverzweigte abelsche Grad-$p$-Erweiterungen beschr\"anken.
\vspace{4mm} \item Im folgenden sei $A(K)=Cl(K)_p$ und $H(K)=H(K)_p$ f\"ur eine feste Primzahl $p$.
\end{itemize}
\end{frame}

\begin{frame}[shrink]
\begin{center}
\Large {\bf Chevalleys Theorem } 
\end{center}
\pause
\begin{theorem}[Furtw\"angler, Chevalley] Sei $L/K$ eine zyklische Erweiterung, in der h\"ochstens ein Primideal verzweigt. Sei $G=Gal(L/K)$ erzeugt von $\sigma\in{G}$, $s=\sigma-1$. Dann gilt 
\begin{eqnarray} KerN_{L/K} &=& A(L)^s,
\end{eqnarray}
und damit: $rk(A(L))\le [L:K]\cdot rk(A(K))$.  
\end{theorem}
\pause
\begin{itemize} [<+->]
\item Haben u.a. erstmals gezeigt, dass $A(L)^G/\imath_{L/K}(A(K))\cong \mathcal{H}^0(G,\mathcal{O}_L^{*})$, falls $L/K$ unverzweigt ist.
\item Dies hat zu einem besonders kurzen und modernen Beweis von Chevalleys Theorem gef\"uhrt.
\item Wir sagen: Erweiterungen, die (1) erf\"ullen, besitzen die Furtw\"angler-Eigenschaft, kurz $F$-Eigenschaft.
\end{itemize}
\end{frame}

\begin{frame}[shrink]
\begin{theorem} Sei $L/K$ wie oben und vom Grad $p$. Dann existiert ein System $\{b_1,...,b_n\}$ in $A(L)$ mit $\mathbb{Z}[s]$-Zyklen $B_i=b_i^{\mathbb{Z}[s]}$, so dass
\begin{eqnarray} KerN_{L/K}=B_1^s\times ...\times B_n^s.
\end{eqnarray}
\end{theorem}
\pause
\begin{definition} Ein $\mathbb{Z}[s]$-Zyklus $B_i$ wie oben hei\ss t {\bf exakt}, falls $B_i\cap A(L)^s=B_i^{s}$ und {\bf nicht exakt} andernfalls.
\end{definition}
\pause
\begin{itemize}[<+->]
\item In nicht-exakten Zyklen gibt es keine Kapitulation:
\end{itemize} 
\pause
\begin{theorem} Sei $b\in{A(L)}$ mit $N_{L/K}(b)=a\in{A(K)}$ und $B=b^{\mathbb{Z}[s]}$. Falls $B$ nicht exakt ist, folgt: $ord(a)=ord(\imath_{L/K}(a))$.
\end{theorem}
\end{frame}

\begin{frame}[shrink]
\begin{center}
\Large {\bf Wachstum von Idealklassen} 
\end{center}
\pause
\begin{definition} Sei $L/K$ wie oben, $G=Gal(L/K)$ erzeugt durch ein $\sigma\in{G}$ und $s=\sigma-1$. Wir definieren $A(K)'=N_{L/K}(A(L))$. Dann sagen wir $L/K$ hat:
\pause
\vspace{3mm} \\ (a) {\bf stabiles Wachstum}, falls $rk(A(K)')=rk(A(L))$ und $rk(A(K))=rk(A(K)^p)$;
\pause \vspace{2mm} \\ \noindent
(b) {\bf zahmes Wachstum}, falls ...
\pause \vspace{2mm} \\ \noindent
(c) {\bf semi-stabiles Wachstum}, falls $A(L)^{s^{p-1}}=\{1\}$;
\pause \vspace{2mm} \\ \noindent
(d) {\bf wildes Wachstum},  sonst.
\end{definition}
\end{frame}

\begin{frame}[shrink]
\pause 
\begin{itemize}[<+->]
\item MAGMA: Alle vier Typen von Wachstum treten auf.
\vspace{2mm}\item In den F\"allen (a)-(c) gilt: $exp(KerN_{L/K})\le p$.
\vspace{2mm}\item Wildes Wachstum: Haben explizit Familien von Erweiterungen konstruiert, in denen $exp(KerN_{L/K})$ beliebig gro\ss \ ist.
\end{itemize} 
\pause

\begin{theorem} In der Situation wie oben sei $b\in{A(L)}$ mit $N_{L/K}(b)=a$. Falls $a$ zu einem minimalen Erzeugendensystem von $A(K)'$ verl\"angerbar ist, dann gilt in den F\"allen (a)-(c):
\begin{eqnarray} ord(b) &=& p\cdot ord(\imath_{L/K}(N_{L/K}(b))).
\end{eqnarray}
\end{theorem}
\pause
\begin{itemize} [<+->]
\item Es folgt in diesem Fall, dass $a$ genau dann in $L$ kapituliert, falls $ord(b)=ord(a)$.
\end{itemize}
\end{frame}

\begin{frame}[shrink]
\begin{center}
\Large {\bf Struktur des Kapitulationskerns} 
\end{center}

\pause
\begin{beispiel} Sei $K=\mathbb{Q}(\alpha)$ mit $\alpha^2=3299$. Dann liefert MAGMA:
\pause 
\begin{itemize} [<+->]
\item $A(K)=<a_1,a_2>\cong C_3\times C_9$
\item Sei $L=H(K)^{<a_1>}$. Dann ist $L/K$ zyklisch vom Grad 9.
\item $P_K(L)=<a_1,a_2^3>\cong C_3\times C_3$
\end{itemize}
\end{beispiel}
\pause
\begin{itemize} [<+->]
\item Dieses Beispiel zeigt, dass sich die Galoisgruppe i.A. nicht in den Kapitulationskern einbetten l\"asst.
\item Mit Hilfe von Galois Kohomologie haben wir u.a. hinreichende Bedingungen gegeben, unter denen dies m\"oglich ist. Diese sind allerdings recht technisch und schwer zu verifizieren.
\item Um mehr \"uber die Struktur des Kapitulationskerns zu erfahren, nehmen wir im folgenden eine zus\"atzliche $G$-Wirkung auf der Idealklassengruppe von $K$ an:
\end{itemize}
\end{frame}

\begin{frame}[shrink]
\begin{center}
\Large {\bf $G$-Wirkung auf der Idealklassengruppe} 
\end{center}
\pause
Sei $K/k$ eine abelsche Galoiserweiterung mit Galois Gruppe $G=Gal(K/k)$ und $p\nmid |G|$.
\vspace{4mm}
\pause \begin{itemize} [<+->]
\item Sei $\alpha\in{\mathbb{Z}_p[G]}$ eine Idempotente. 
\item Dies liefert eine Zerlegung von $\mathbb{Z}_p[G]$:
\[\mathbb{Z}_p[G]=\alpha\mathbb{Z}_p[G]\oplus (1-\alpha)\mathbb{Z}_p[G].\]
\item Und eine Zerlegung von $A(K)$:
\[A(K)=A(K)_{\alpha}\times A(K)_{1-\alpha},\ \mbox{wobei}\]
\[A(K)_{\alpha}=A(K)^{\alpha}=\{a^{\alpha}|\ a\in{A(K)}\}\]
\item Wir nennen $A(K)_{\alpha}$ eine irreduzible Komponente von $A(K)$, falls $\alpha$ primitiv ist.
\end{itemize}
\end{frame}

\begin{frame}[shrink]
\vspace{4mm}
\begin{theorem} Sei $K/k$ wie oben. Dann kann $A(K)$ in ein direktes Produkt von zyklischen $\mathbb{Z}_p[G]$-Untermoduln ($\mathbb{Z}_p[G]$-Zyklen) zerlegt werden. 
\end{theorem}
\pause
\vspace{4mm}
\begin{korollar} Sei $B$ ein $\mathbb{Z}_p[G]$-Zyklus in einer irreduziblen Komponente $A(K)_{\alpha}$. Dann:
\pause 
\begin{itemize} [<+->]
\vspace{2mm}\item $\mathbb{Z}_p[G]$ wirkt transitiv auf Elemente in $B$ von gleicher Ordnung.
\vspace{2mm}\item Sei $L\supset K\supset k$, $L/k$ Galois, $n\in{\N}$. Dann kapituliert entweder jede Idealklasse $a\in{B}$ mit $ord(a)=p^n$ in $L$ oder keine.
\vspace{2mm} \item Unter gewissen Voraussetzungen kapitulieren alle $a\in{A(K)}$ mit $ord(a)\le p^n$ in einer Teilerweiterung von $H(K)/K$ vom Exponent $p^n$.
\vspace{2mm}\item $KH(k)$ ist gerade der Genusk\"orper von $K/k$.
\end{itemize}
\end{korollar} 
\end{frame}

\begin{frame}[fragile,allowframebreaks]
\begin{center}
\Large {\bf Kapitulation in Imagin\"ar Quadratischen K\"orpern}
\end{center}

\begin{longtable}{|c|c|c|} \hline  Nr. & Discriminant & $(P_K(L_1),...,P_K(L_6))$\\  \hline\hline  1 & -12451 & (3,6,4,1,5,2)\\ \hline 2 & -17944 & (2,1,5,4,6,3)\\ \hline 3 & -30263 & (6,5,3,4,2,1)\\ \hline 4 & -33531 & (3,2,4,6,1,5)\\ \hline 5 & -37363 & (1,2,3,4,6,5)\\ \hline 6 & -38047 & (3,1,6,4,2,5)\\ \hline 7 & -39947 & (5,1,4,6,2,3)\\ \hline 8 & -42871 & 
(2,1,6,5,4,3)\\ \hline 9 & -53079 & (2,6,4,1,5,3)\\ \hline 10 & -54211 & (2,6,3,4,1,5)\\ \hline 11 & -58424 & (3,2,1,5,4,6)\\ \hline 12 & -61556 & (6,5,4,1,3,2)\\ \hline 13 & -62632 & (2,6,6,6,6,6)\\ \hline 14 & -63411 & (5,2,6,4,1,3)\\ \hline 15 & -64103 & (1,3,6,5,2,4)\\ \hline 16 & -65784 & (3,6,4,5,2,1)\\ \hline 17 & -66328 & (6,2,4,3,5,1)\\ \hline 18 & -67031 & (4,4,4,1,4,4)\\ \hline 19 & -67063 & (1,1,1,6,1,1)\\ \hline 20 & -67128 & (3,2,4,6,1,5)\\ \hline 21 & -69811 & (3,1,4,2,5,6)\\ \hline 22 & -72084 & (3,1,5,6,4,2)\\ \hline 23 & -74051 & (2,1,3,6,5,4)\\ \hline 24 & -75688 & (5,6,5,5,5,5)\\ \hline 25 & -81287 & (4,1,3,5,2,6)\\ \hline 26 & -83767 & (5,4,3,2,1,6)\\ \hline 27 & -84271 & (5,4,2,1,6,3)\\ \hline 28 & -85099 & (6,2,1,4,3,5)\\ \hline 
\end{longtable}
\end{frame}

\begin{frame}[shrink]
\vspace{4mm}
\begin{theorem} Sei $K$ ein imagin\"ar quadratischer K\"orper mit $rk(A(K))=2$, $p>3$. Seien $L_1,...,L_{p+1}$ die Zwischenk\"orper von $H(K)/K$ vom Grad $p$ \"uber $K$, so dass $rk(KerN_{L_i/K})\le rk(KerN_{L_j/K})$, $\forall\ i<j$. Wir nehmen an: 
\pause \vspace{2mm} \\ \noindent (A1) $exp(KerN_{L_j/K})=p$, $\forall\ j$;
\pause \\ (A2) $rk(KerN_{L_1/K})=2$.
\vspace{2mm} \\ 
\pause \noindent Dann: (i) $P_K(L_i)\not=P_K(L_j)$, $\forall\ i\not=j$ $\Leftrightarrow\  rk(KerN_{L_u/K})=2$, $\forall\ u$.
\pause \vspace{2mm} \\ \noindent (ii) $P_K(L_i)=P_K(L_j)$, $\forall\ 1\le i\not=j\le p$, andernfalls.
\end{theorem}
\pause
\begin{itemize} \item Im ersten Fall sagen wir $K$ hat {\bf 1-1-Kapitulation} und im zweiten Fall sagen wir $K$ hat {\bf $p$-Kapitulation}.
\end{itemize}
\end{frame}

\begin{frame}[shrink]
\vspace{6mm}
\begin{itemize} [<+->]
\item Wir haben gezeigt, dass $rk(KerN_{L_j/K})>0$ gerade ist f\"ur alle $j$.
\vspace{6mm}\item Die Wahrscheinlichkeit, dass die Annahmen (A1) und (A2) erf\"ullt sind, sollte f\"ur wachsendes $p$ gegen 1 gehen. 
\vspace{6mm}\item Um dies zu untermauern, haben wir eine Heuristik \"uber die Klassengruppe von unverzweigten zyklischen Erweiterungen von imagin\"ar quadratischen Zahlk\"orpern entwickelt.
\vspace{6mm}\item Das obige Theorem l\"asst sich in gewisser Weise f\"ur $rk(A(K))>2$ verallgemeinern.
\end{itemize}
\end{frame}

\begin{frame}[shrink]
\begin{center}
\large \textbf{Heuristiken \"uber die Klassengruppe von Unverzweigten Zyklischen Erweiterungen von Imagin\"ar Quadratischen Zahlk\"orpern}
\end{center}
\pause
Sei $L/K$ unverweigte zyklische Erweiterung vom Grad $p$, $K$ imagin\"ar quadratisch. Der Einfachheit wegen sei $rk(A(K))=2$. 
\pause
\begin{itemize} [<+->]
\item Frage: Wie h\"aufig tritt der Fall $rk(KerN_{L/K})=2k$ auf, wenn $K$ die imagin. quadr. K\"orper mit $rk(A(K))=2$ durchl\"auft?
\vspace{4mm}
\item Heuristik st\"utzt sich auf Ideen von Cohen-Lenstra unter Ber\"ucksichtigung der entwickelten Resultate ($KerN_{L/K}$ hat geraden Rang und ist nach Chevalley $\Z[s]$-zyklisch).
\vspace{4mm}
\item Da wir nur an $KerN_{L/K}$ interessiert sind, ist nur der Rang von $A(K)$ entscheidend und nicht etwa der Exponent. 
\vspace{4mm} \item Heuristik: F\"ur $p\rightarrow \infty$, geht die rel. H\"aufigkeit f\"ur $rk(KerN_{L/K})=2$ gegen 1.
\end{itemize}
\end{frame}

\begin{frame}[shrink]
\begin{center}
\large \textbf{Heuristik und Numerische Daten im Vergleich}
\end{center}
\pause
\begin{itemize} [<+->]
\item Heuristik ist in guter \"Ubereinstimmung mit den numerischen Daten: Im Fall $p=3$ haben wir 12484 solcher Erweiterungen $L/K$ untersucht und mit MAGMA folgende relative H\"aufigkeiten (mittlere Spalte) ermittelt:
\end{itemize}
\vspace{4mm}
\pause
\begin{longtable}{|c|c|c|} \hline  & $rk(KerN_{L/K})=2k$ & Heuristik \\  \hline\hline  $k=1$ & 0.8992 & 0.8889\\ \hline  $k=2$ & 0.0950 & 0.0989\\ \hline $k=3$ & 0.0071 & 0.0110\\ \hline
\end{longtable}
\end{frame}

\end{document}